\documentclass[runningheads,a4paper]{llncs}

\usepackage{amssymb}
\usepackage{amsmath}
\usepackage{graphicx}
\usepackage{comment}
\usepackage{algorithmic}
\usepackage{algorithm} 
\usepackage{setspace}
\usepackage{listings}

\usepackage{color} 
\definecolor{mygreen}{RGB}{28,172,0} 
\definecolor{mylilas}{RGB}{170,55,241}

\usepackage{url}
\usepackage{hyperref}

\usepackage{tikz}
\usetikzlibrary{external}
\urldef{\mailsa}\path|lmarcin@mimuw.edu.pl|
\urldef{\mailsb}\path|jan.valdman@utia.cas.cz |

\newcommand{\nb}{n_b}
\newcommand{\dx}{\,\mathrm{d}x}

\newcommand{\Mod}[1]{\ (\mathrm{mod}\ #1)}

\begin{document}

\mainmatter  
\title{MATLAB Implementation of Element-based Solvers}

\titlerunning{MATLAB Implementation of element-based solvers}

\author{Leszek Marcinkowski$^1$
\thanks{The work of the 1st author was partially supported by Polish Scientific Grant: National Science Center  2016/21/B/ST1/00350.}
and Jan Valdman$^2$
\thanks{The work of the 2nd and the corresponding author was supported by the Czech Science Foundation (GACR), through the grant 
GA17-04301S.}
}

\authorrunning{Leszek Marcinkowski and Jan Valdman}

\institute{$^1$ Faculty of Mathematics, University of Warsaw, 
Warszawa, Poland\\
\mailsa\\ 
$^2$
Faculty of Science, University of South Bohemia, 
\v{C}esk\'{e}~Bud\v{e}jovice, \\ and The Czech Academy of Sciences, Institute of Information Theory and Automation, 
Prague, Czechia \\
\mailsb
}
\maketitle
\begin{abstract}
Rahman and Valdman (2013) introduced a vectorized way to assemble finite element stiffness and mass matrices in MATLAB. Local element matrices are computed all at once by array operations and stored in multi-dimentional arrays (matrices). We build some iterative solvers on available multi-dimentional structures completely avoiding the use of a sparse matrix. 
\keywords{MATLAB code vectorization, finite elements, stiffness and mass matrices, iterative solvers}
\end{abstract}

\section{Motivation example}
\begin{figure}[t]
  \centering
  \includegraphics[width=0.7\textwidth]{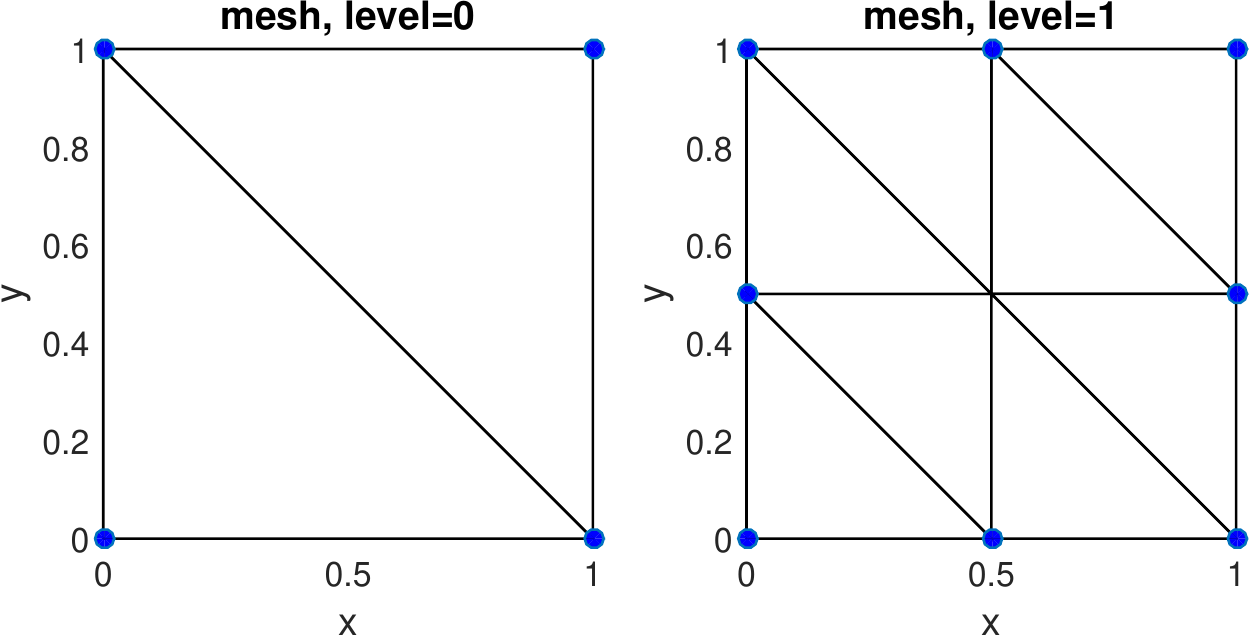}
  \caption{Two examples of triangular meshes of a unit square domain $\Omega$ with $n_e=2$ elements and $n_n=4$ nodes (left) and $n_e=8$ elements and $n_n=9$ nodes (right).}
  \label{fig:k}
\end{figure}

We solve a benchmark boundary value problem
$$ -\triangle u + \nu  u = f \qquad \mbox{on } x \in \Omega = (0,1) \times (0,1)$$ 
for given $f \in L^2(\Omega)$ and a parameter $\nu \geq 0$. Nonhomogeneous Dirichlet or homogeneous Neumann boundary conditions are assumed on parts of boundary $\partial \Omega$  and measure of the Dirichlet boundary has to be positive for $\nu=0$. 
A finite element method is applied and leads to a linear system of equations
\begin{equation}
A u = (K + \nu M ) u = b \label{linear_system}
\end{equation}
for an unknown vector $u \in \mathbb{R}^{n_n}$, where $n_n$ denotes the number of mesh nodes (vertices). Stiffness and mass matrices $K, M \in \mathbb{R}^{n_n \times n_n}$ and the right hand side vector $b \in \mathbb{R}^{n_n}$ are defined as
\begin{equation} \label{matrices_K_M}
K_{ij} =\int_\Omega \nabla \Phi_i \cdot \nabla \Phi_j  \dx, \qquad 
M_{ij} =\int_\Omega \Phi_i \, \Phi_j \dx, \qquad 
b_{j} =\int_\Omega f \, \Phi_j \dx
\end{equation}
using local basis functions $\Phi_i$  for $i=1,\dots, n_n$ and $\nabla$ denotes the gradient operator. Fig. \ref{fig:k} shows an example of a 2D discretization of $\Omega$. 

Sparse matrices $K, M$ are generated as 
\begin{equation} 
K = \sum_{e=1}^{n_e} C^T_e K_e C_e, \qquad 
M = \sum_{e=1}^{n_e} C^T_e M_e C_e,
\end{equation}
where $n_e$ denotes a number of mesh elements (number of triangles in Fig.  \ref{fig:k}), 
$$K_e, M_e \in \mathbb{R}^{\nb \times \nb}, \qquad e=1,\dots, n_e$$ are local element matrices 
and 
$$C_e \in \mathbb{R}^{\nb \times n_n}, \qquad e=1,\dots, n_e$$ 
are Boolean connectivity matrices which distribute the terms in local element matrices to their associated global degrees of freedom. Here, $\nb$ denotes a number of local basic functions. In the simplest case of nodal linear ($P_1$) finite elements: 
\begin{eqnarray*}
&&\nb=3 \qquad \mbox{for triangles in 2D},\\
&&\nb=4 \qquad \mbox{for tetrahedra in 3D}.
\end{eqnarray*}
Extensions to higher order isoparametric elements are also possible. 
\begin{figure}
  \centering
  \includegraphics[width=0.7\textwidth]{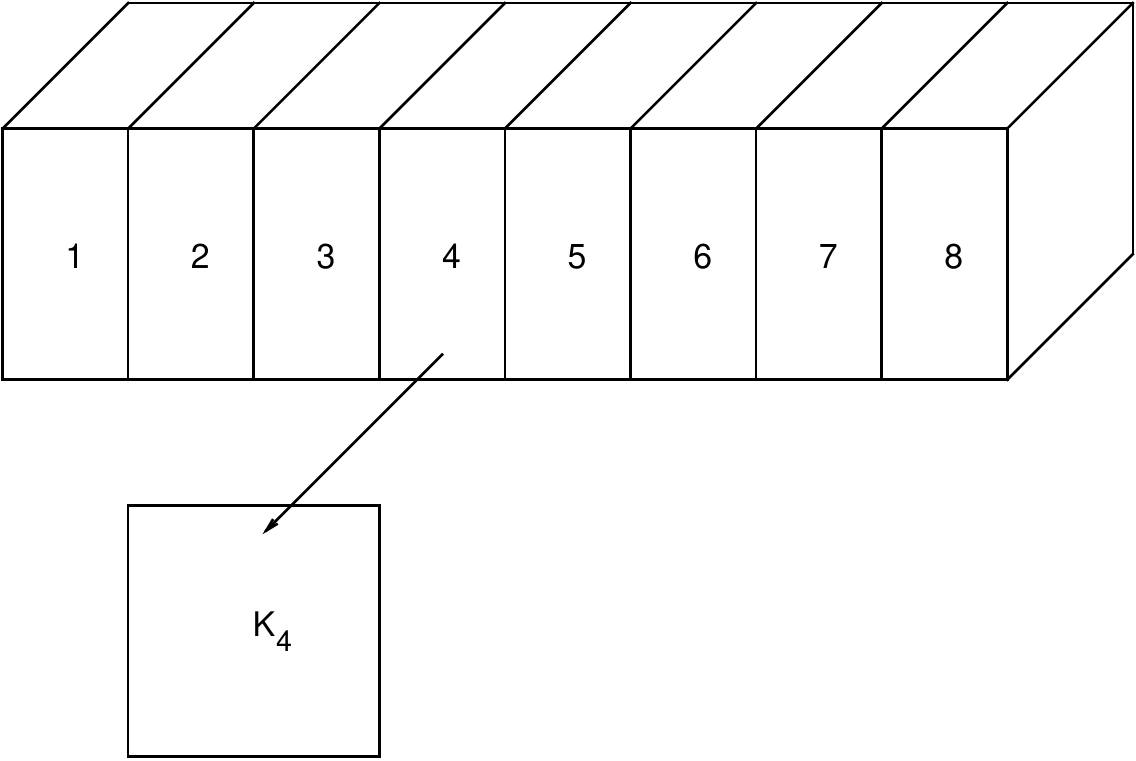}
  \caption{Example of a 3-dimensional array storing all local stiffness matrices. The matrix corresponds to a triangular mesh with 8 elements displayed on Figure \ref{fig:k} (right). A particular local stiffness matrix $K_4 \in \mathbb{R}^{3 \times 3}$ is indicated. }
  \label{fig:3Dmatrix} 
\end{figure}
All matrices $K_e, M_e$ for  $e=1,\dots, n_e$ are generated at once using 
vectorized routines of \cite{RahmanValdman2013}. 
They are stored as 3-dimensional full matrices (see Figure \ref{fig:3Dmatrix}) of sizes
$$ \nb \times  \nb \times n_e.$$ 
The storage of 3-dimensional matrices contains certain memory overheads in comparison to sparse matrices (which can be automatically generated from them), since local contributions from restrictions of basis functions to shared elements are stored separately. 
Our aim is to build and explain in detail simple linear iterative solvers based on local element matrices $K_e, M_e$ without assembling the sparse matrices $M, K$. This is our first attempt in this direction and therefore we show the possibility of this approach rather than efficient implementations and runtimes. The complementary software to this paper is available for download 
\begin{center}
\url{https://www.mathworks.com/matlabcentral/fileexchange/70255} .
\end{center}

\section{Element-based solvers}
Some examples of element-based iterative solvers are provided including their simple MATLAB implementations. All are based on a vectorized computation of a (column) residual vector
\begin{equation} \label{residual}
r := b - A x 
\end{equation} 
for a given approximation (column) vector $x \in \mathbb{R}^{n_n}$.
\begin{algorithm}
\caption{residual computation - looped version}
\label{residduum_update_elementwise}
\hspace{0.2cm}
\begin{spacing}{1.2}
\begin{algorithmic}[1]
  \FOR{$e=1,\dots, n_e$} 
    \STATE $x_e = R_e x, \qquad \qquad  \qquad \qquad $ (restriction)
     \STATE $r_e = b_e - A_e x_e, \qquad \qquad  \quad\enspace \, $ (local residual)
  \ENDFOR 
   \STATE $r = \sum_{e=1}^{n_e} C^T_e r_e. \qquad \qquad  \quad \enspace \enspace \enspace $ (assembly)
\end{algorithmic}
\end{spacing}
\end{algorithm}
The residual is computed using local matrices and local vectors
$$A_e := K_e + \nu M_e \in \mathbb{R}^{\nb \times \nb}, \quad b_e \in \mathbb{R}^{\nb}, \qquad e=1,\dots, n_e.$$
Matrices $R_e \in \mathbb{R}^{\nb \times n_n}, e=1,\dots, n_e$ 
are restriction matrices from global to local indices.
 Note that elementwise  evaluations inside the loop (lines 2 and 3) operate with local matrices and local vectors only. A fully vectorized MATLAB version of Algorithm \ref{residduum_update_elementwise} follows:
\begin{lstlisting}
function r=residual_e(A_e,bt_e,x,ind_e,indt)
x_e=x(ind_e);                               %restriction - all  
rt_e=bt_e-avtam(x_e,A_e);                   %residual - all
r=accumarray(indt(:),rt_e(:));              %assembly - all
end     
\end{lstlisting}
Clearly, matrices $R_e$ and $C_e$ of Algorithm \ref{residduum_update_elementwise} 
are not stored, but their operations are replaced by a convenient indexing using two index arrays: 
$$\verb+ind_e+ \in \mathbb{I}^{\nb \times 1 \times n_e}, \quad \verb+indt+ \in \mathbb{I}^{\nb \times n_e}.$$
Both arrays contain the same global nodes numbers corresponding to each element, but they are ordered differently with respect to their operations.

All objects indexed by elements are stored as full higher dimensional matrices and their names end with a symbol $\verb+_e+$. 

\subsection{Richardson iteration}
We recall few examples of iterative methods based on a residual computation 
More details about them can be found eg. in \cite{Hackbush2014,SamarskiiNikolaev1989}.
One of the simplest iterative methods to solve 
\eqref{linear_system} is the Richardson iteration for iterations $k=0, 1,2,\dots$ in the form
\begin{equation} \label{richardson}
\begin{split}
    &r^{k}=b - A x^{k}, \\
    &x^{k+1} = x^{k} + \omega \, r^{k} 
\end{split}
\end{equation} 
with the initial column vector $x^0 \in \mathbb{R}^n$ and a given positive parameter $\omega > 0$. 
The optimal coefficient is equal to 
$\omega_{opt}=\frac{1}{\lambda_1+\lambda_2}$ for $A=A^T>0$, where $\lambda_1$ is the smallest and $\lambda_2$ the largest eigenvalue of $A$. 
For this $\omega_{opt}$ the convergence estimate is the fastest, i.e. 
$$\|x^k-u\|_2\leq \frac{\lambda_2-\lambda_1}{\lambda_2+\lambda_1}\|x^{k-1}-u\|_2.$$ Here, $u \in \mathbb{R}^n $ is the solution of \eqref{linear_system}. A MATLAB version follows:
\begin{lstlisting}
function x=Richardson_e(A_e,bt_e,x0,ind_e,indt,iters,lam2,lam1,nd)
omega=2/(lam2+lam1);                        %optimal parameter
x=x0;                                       %iteration initial
for k=0:iters-1
    r=residual_e(A_e,bt_e,x,ind_e,indt);    %residual comput. 
    r(nd)=0;                                %dirichlet condit.
    x=x+omega*r;                            %iteration update  
end
end
\end{lstlisting}

\subsection{Chebyshev iteration}
The Chebyshev polynomial (of the first kind) of degree $N \in \mathbb{N}_0$  is defined by
$$
 T_N(x):=\cos(N\arccos (x)), \qquad x\in [-1,1]
$$
and it is known to have roots in the form
$$\alpha_k=\cos (\pi (k+1/2)/N), \qquad  k=0,\ldots,N-1.$$
Consequently, a shifted and scaled polynomial 
$$
P_{N}(t)=T_{N}\left( \left(\frac{-2}{\lambda_2-\lambda_1}\right)\left(t-\frac{\lambda_1+\lambda_2}{2}\right)\right)/{C_N}, \qquad t\in [\lambda_1,\lambda_2],
$$
with the scaling factor  $C_N=T_{N}\left(\frac{\lambda_1+\lambda_2}{\lambda_2-\lambda_1}\right)$
satisfies the condition $P_{N}(0)=1$. It also has $N$ distinct roots 
$$\alpha_k=\frac{\lambda_1+\lambda_2}{2} - \frac{\lambda_2-\lambda_1}{2}\cos \left(\frac{\pi (k+1/2)}{N}\right), \qquad k=0,\ldots,N-1$$ 
lying in $(\lambda_1,\lambda_2)$. This polynomial has the smallest maximum norm on $[\lambda_1,\lambda_2]$ over all polynomials of degree less or equal $N$ which are equal to one at zero. 
\subsubsection{Two-level Chebyshev iteration}
The cyclic two-level Chebyshev  iterative methods to solve 
\eqref{linear_system} is in the form
\begin{equation} \label{Czebyshev2lev}
\begin{split}
    &r^{k}=b - A x^{k}, \\
    &x^{k+1} = x^{k} + \alpha_{k \Mod N}^{-1} r^{k}.
\end{split}
\end{equation} 
The method is convergent if all eigenvalues of $A$ are contained in $[\lambda_1,\lambda_2] \subset (0,\infty)$. The optimal convergence is accessed  where  $\lambda_1$ is the minimal eigenvalue and $\lambda_2$ the maximal eigenvalue of $A$.
Note that after $N$ iterations we get 
\begin{equation} \label{Cheberror1}
    x^{N}-u=\Pi_{k=0}^{N-1} (I-\alpha_k^{-1} A)(x^0-u)=P_{N}(A)(x^0-u).
\end{equation}  
and then after $\ell N$ iterations we get $x^{\ell N}-u=(P_{N}(A))^\ell(x^0-u)$.
Note that the Richardson iteration \eqref{richardson} is the special case of this method with $N=1$.
This error formula gives us,
$$
 \|x^{N}-u\|_2\leq \sum_{t \in [\lambda_1,\lambda_2]} |P_{N}(t)|\|x^0-u\|_2 \leq 2\left(\frac{\sqrt{\lambda_2}-\sqrt{\lambda_1}}{\sqrt{\lambda_2}+\sqrt{\lambda_1}}\right)^N\|x^0-u\|_2.
$$
A MATLAB version follows:

\begin{lstlisting}
function x=Chebyshev2Level_e(A_e,bt_e,x0,ind_e,indt,iters,lam2,lam1,N,nd)
d=(lam2+lam1)/2; c=(lam2-lam1)/2;
k=0:N-1; alphas=d+c*cos(pi*(1/2+k)/N);
x=x0;                                       %iteration initial
for k=0:iters-1
    r=residual_e(A_e,bt_e,x,ind_e,indt);    %residual comput. 
    r(nd)=0;                                %dirichlet condit.
    alpha=alphas(mod(k,N)+1);
    x=x+(1/alpha)*r;                        %iteration update  
end
end
\end{lstlisting}

\subsubsection{Three-level Chebyshev iteration}
We now present the three-level Chebyshev iteration, cf.  
e.g. \cite{SamarskiiNikolaev1989,Hackbush2014},
The method is defined by the error equation, cf. also (\ref{Cheberror1}),
\begin{equation} \label{Cheberror2}
 x^k-u=P_k(A)(x^0-u), \qquad k=0,1,2,\ldots
\end{equation}
and
its implementation is  based on the following recurrence relation
$$
 T_k(t)=2 t \, T_{k-1}(t)-T_{k-2}(t), \qquad k>1, \qquad T_1(t)=t,\;T_0(t)=1.
$$
This relation for $t_k$ yields the recurrence formula for $k>1$,
\begin{eqnarray}
 P_{k+1}(x)&=&2\frac{\lambda_1+\lambda_2-2x}{\lambda_2-\lambda_1}\frac{C_k}{C_{k+1}}P_k(x) - \frac{C_{k-1}}{C_{k+1}}P_{k-1}(x), \nonumber \\
  C_{k+1}&=&2\frac{\lambda_1+\lambda_2}{\lambda_2-\lambda_1}C_k - C_{k-1}, \label{recrCk}
\end{eqnarray}
where
$$
P_1(x)=C_1^{-1}\frac{\lambda_1+\lambda_2-2x}{\lambda_2-\lambda_1}, \qquad P_0=1, \qquad
 C_1=\frac{\lambda_1+\lambda_2}{\lambda_2-\lambda_1}, \qquad C_0=1. 
$$
For $k=1$ we get $x^1-u=P_1(A)(x^0-u)$ and 
$$
 x^1=u+\frac{\lambda_2-\lambda_1}{\lambda_1+\lambda_2} \frac{\lambda_1+\lambda_2-2A}{\lambda_2-\lambda_1}(x^0-u)=x^0+\frac{2}{\lambda_2-\lambda_1}r_0,
$$
where  $r_0=b-Ax^0$. Note that $x^1$ is computed as one iteration of the Richardson method applied to $x^0$ with the optimal coefficient, cf. (\ref{richardson}). 
Our method is defined by (\ref{Cheberror2}), thus using the above recurrence relation we get for $ k>1$,
$$
 x^{k+1}-u=\frac{2C_k}{C_{k+1}}\left(\frac{\lambda_1+\lambda_2}{\lambda_2-\lambda_1}I -\frac{2}{\lambda_2-\lambda_1}A \right)(x^k-u)
 - \frac{C_{k-1}}{C_{k+1}}(x^{k-1}-u).
$$
Since $$1=2\frac{\lambda_1+\lambda_2}{\lambda_2-\lambda_1}\frac{C_k}{C_{k+1}} - \frac{C_{k-1}}{C_{k+1}}$$ we see that
$$
  x^{k+1}=\frac{2C_k}{C_{k+1}}\frac{\lambda_1+\lambda_2}{\lambda_2-\lambda_1}x^k +\frac{4}{\lambda_2-\lambda_1}\frac{C_k}{C_{k+1}}(b-Ax^k)
 - \frac{C_{k-1}}{C_{k+1}}x^{k-1}
$$
and utilizing this identity once more we have the three level Chebyshev iterations
\begin{eqnarray}\label{Czebyshev3lev}
 x^{k+1}&=&x^k + \frac{C_{k-1}}{C_{k+1}}(x^k-x^{k-1})+\frac{4}{\lambda_2-\lambda_1}\frac{C_k}{C_{k+1}}r_k, \qquad k>1,\\
 x^1&=&x^0+\frac{2}{\lambda_2-\lambda_1}r_0 \nonumber
\end{eqnarray}
with $r_k=b-Ax^k$ $k=0,1,2,\ldots$. We remind that  the scaling factors $C_k$ are defined by (\ref{recrCk}).
Note that $x^N$ in the both 2-level and 3-level iterations,  cf. (\ref{Czebyshev2lev}) and  (\ref{Czebyshev3lev}), are equal to each other 
what follows from \eqref{Cheberror1} and \eqref{Cheberror2}. A MATLAB version reads:

\begin{lstlisting}
function x=Chebyshev3Level_e(A_e,bt_e,x0,ind_e,indt,iters,lam2,lam1,nd)  
d=(lam2+lam1)/2; c=(lam2-lam1)/2;
x=x0;
r=residual_e(A_e,bt_e,x,ind_e,indt);        %residuum comput.
r(nd)=0;                                    %dirichlet condit.
for k = 0:iters-1  
    z=r;
    if (k==0)  
        p=z; alpha=1/d;
    else
        beta=(c*alpha/2)^2;
        p=z+beta*p; alpha=1/(d - beta/alpha);
    end
    x=x+alpha*p;
    r=residual_e(A_e,bt_e,x,ind_e,indt);    %residuum comput.
    r(nd)=0;                                %dirichlet condit.
end
end
\end{lstlisting}

\section{Numerical experiments}
We consider for simplicity the case of the square domain $\Omega = (0,1) \times (0,1),$ no mass matrix ($\nu=0$) and nonhomogenous Dirichlet boundary conditions
$u=1$ for $x \in \partial \Omega.$
    For a uniformly refined triangular mesh (see Figure  \ref{fig:k}) with $n^2$ nodes (also counting boundary nodes), there are $(n-2)^2$ eigenvalues of $A=K$ in the form
$$ \lambda = 4 \left(\sin^2{\frac{i \pi}{2(n-1)}} + \sin^2{\frac{j \pi}{2(n-1)}}  \right), \quad i,j=1,\dots,n-2$$
and the minimal eigenvalue $\lambda_1$ is obtained for $i=j=1$ and the maximal eigenvalue $\lambda_2$ for $i=j=n-2$. We utilize these eigenvalue bounds for all mentioned iterations methods. Furthermore, we assume a constant function $f=1$ for $x \in \Omega$. 
\begin{figure}[h]
  \centering
  \includegraphics[width=\textwidth]{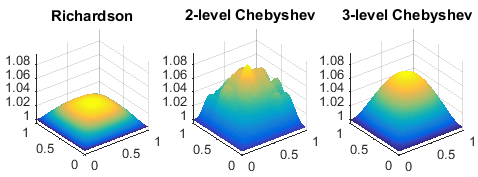}
  \caption{Final iterates.}
  \label{fig:iterations}
\end{figure}

For a given number of iterations (we choose 124 iterations) and a mesh with $1089=33^2$ nodes, final iterates are displayed in Figure \ref{fig:iterations}. Only the 3-level Chebyshev method converged optically to the exact solution. Richardson requires more steps to improve its convergence and the 2-level Chebyshev (with $N=32$) demonstrates a known instability. The remedy to fix this instability would be to reorder values of precomputed parameters $\alpha_k$ to enhance the stability. Time performance is also reasonable for finer meshes. On a  (level 10) mesh with 
$1050625=1025^2$ nodes, the assemblies of 3-dimensional arrays $K_e, M_e$ take around 5 seconds each and 124 iterations take around 50 seconds for all iteration methods. The direct solver of MATLAB takes 5 seconds. Since number of iterations to obtain a convergence with respect to a given tolerance is known to grow as a function of condition number of finer meshes, we need to combine studied solvers with preconditioners or use several iterations of them as smoothers for instance in multigrid procedures. 

\section*{Outlooks}
We are interested in developing preconditioners for discussed solvers 
on multi-dimensional structures and extension to edge elements based on \cite{AnjamValdman2015}.

\bibliographystyle{abbrv}

\end{document}